\documentclass[preprint,1p]{elsarticle}
\usepackage{amsmath}
\usepackage{amssymb}

\usepackage[export]{adjustbox}
\usepackage[usenames,dvipsnames]{color}

\newcommand{\firstRev}[1]{\color{black}{#1}\color{black}}
\newcommand{\secondRev}[1]{\color{black}{#1}\color{black}}

\journal{Computers and Mathematics with Applications}
\begin{document}
\begin{frontmatter}

\title{Tracer Transport within an Unstructured Grid Ocean Model using Characteristic Discontinuous Galerkin Advection}
\author[CCS]{D.~Lee\corref{cor}}
\ead{davelee2804@gmail.com}
\author[CCS]{M.~Petersen}
\author[CCS]{R.~Lowrie}
\author[T]{T.~Ringler}

\address[CCS]{Computer, Computational and Statistical Sciences, Los Alamos National Laboratory, Los Alamos, NM 87545, USA}
\address[T]{Theoretical Division, Los Alamos National Laboratory, Los Alamos, NM 87545, USA}
\cortext[cor]{Corresponding author. Tel. +61 452 262 804.}

\begin{abstract}

In a previous article a characteristic discontinuous Galerkin (CDG) advection scheme was presented for tracer transport \cite{Lee16}. 
The scheme is conservative, unconditionally stable with respect to time step and scales sub-linearly with the number of 
tracers being advected. Here we present the implementation of the CDG advection scheme for tracer transport within 
MPAS-Ocean, a Boussinesque unstructured grid ocean model with an arbitrary Lagrangian Eulerian vertical coordinate. The 
scheme is implemented in both the vertical and horizontal dimensions, and special care is taken to ensure that the
scheme remains conservative in the context of moving vertical layers. Consistency is ensured with respect to the dynamics
by a renormalization of the fluxes with respect to the volume fluxes derived from the continuity equation.
For spherical implementations, the intersection of the flux swept regions and the Eulerian grid are determined for 
great circle arcs, and the fluxes and element assembly are performed on the plane via a length preserving projection.
Solutions are presented for a suite of test cases and comparisons made to the existing flux corrected transport scheme in 
MPAS-Ocean.

\end{abstract}


\end{frontmatter}

\section{Introduction}

Tracer advection constitutes a large portion of the compute time for modern global climate models, due to 
the large number of chemical and hydrometeor species that must be accounted for. For physical consistency, 
the advection of tracers must be conservative while being as numerically accurate and computationally 
efficient as possible. In a previous article, a novel characteristic discontinuous Galkerin (CDG) advection 
scheme was presented which allows for conservative advection on unstructured meshes and arbitrarily long 
time steps while also scaling sub-linearly with the advection of additional tracers \cite{Lee16}. In 
idealized geometries, the CDG scheme was found to outperform a traditional flux corrected transport (FCT) 
scheme \cite{BB73, Zalesak79} for a moderate number of tracers, and was shown to converge at higher order, 
through the use of a modal basis expansion of the tracer in each element.

Several conservative transport schemes have recently been presented based on the idea of integrating either
the edges or the entire element backwards along Lagrangian characteristics and integrating over the resultant 
area in order to determine either the fluxes or the element values at the new time level. These schemes have
the appealing properties that they are unconditionally stable with respect to time step, and so may be run
with longer time steps than the underlying dynamics requires, and their performance scales sub-linearly with 
the number of tracers being advected, due to the fact that the computation of the Lagrangian characteristics 
may be reused for additional tracers. The Incremental Remap (IR) scheme, which has been implemented on both 
planar Cartesian \cite{DB00} and spherical geodesic \cite{LR05} grids, exploits this idea in flux form. The 
IR scheme uses the mean values of the tracer in the neighbouring elements in order to reconstruct the tracer 
gradients, which are in turn used to integrate the swept regions of the edges. It has previously 
\firstRev{been } shown
\cite{Lee16} that the CDG scheme is approximately as accurate as the IR scheme at half the resolution, due
to the compact nature of the trial functions compared to the IR gradients.

In remap form, by which the entire element is integrated backwards along the characteristics via the Reynolds
transport theorem in order to compute its new conservative value, this idea has been used to construct the  
Conservative Semi-Lagrangian Multi-tracer transport scheme (CSLAM) \cite{LNU10,ELT13}. In this formulation
the intersection of the pre-image of the element and the Eulerian grid at the previous time level is integrated 
using line integrals via the Gauss-Green theorem in order to determine the weights of a quadratic polynomial 
representation of the tracer in each element. This scheme is currently in use within the High-Order Methods 
Modelling Environment (HOMME) atmospheric model \cite{TF10}. 

The methods discussed above use some form of reconstruction to determine the higher order structure of the 
tracer field. An alternative approach is to introduce a set of test functions which are integrated 
along velocity characteristics so to satisfy the adjoint equation to the weak form of the problem. This is 
the approach used in the Eulerian-Lagrangian Localized Adjoint Method (ELLAM) \cite{CRHE90,RC02}.
One downside of the ELLAM method however is that it requires the assembly and solution of a global system of 
equations, as either a finite volume or finite element problem. This issue is negated in a similar and recently
developed semi-Lagrangian discontinuous Galerkin scheme for tracer transport in atmospheric flows
on cubed spheres \cite{GNQ14}, which prognoses the trial function representation 
of the tracer, while integrating the quadrature points of the test functions forwards in time along
velocity characteristics in order to satisfy the adjoint equation. The method is applied in one dimension, for 
which the quadrature points of the Lagrangian pre-image of the element are integrated forward in time, where 
they are used to evaluate the tracer via its trial function representation. This is required in order to 
preserve the values of the test functions along characteristics, also a feature of the CDG scheme.

The CDG scheme is similar to the previous semi-Lagrangian discontinuous Galerkin scheme \cite{GNQ14}
in that the higher order structure is prognosed via the solution of a system of linear equations
for the coefficients of the trial functions in each element, with the fluxes determined via an integration
of the swept region made by the vertices of the edges along characteristics.
However unlike the previous scheme it may be applied in two dimensions without
the use of dimensional splitting, and so is suitable for fully unstructured grids. 
In this paper we present the implementation of the CDG scheme within the MPAS-Ocean model \cite{Ringler13},
a mimetic c-grid finite volume model for the advection of both passive and active tracers.
The scheme is implemented in both the horizontal, on planar and spherical unstructured grids, and in
the vertical, which makes use of an arbitrary Lagrangian Eulerian (ALE) grid \cite{Petersen15}.
Consistency between the dynamics and the transport scheme is ensured via a normalization of the fluxes
by the volume fluxes derived from the continuity equation, and 
special care is taken to ensure that the vertical advection remains conservative in the context of the
moving layers.

The remainder of this paper is presented as follows: In section 2 the formulation of the CDG scheme for the 
advection equation is presented, with particular emphasis on the 
\firstRev{splitting of fluxes between the horizontal and vertical dimensions and the construction of the 
vertical scheme in the context of the ALE vertical coordinate. }
\firstRev{Section 3 } presents the results of various idealized test cases
and comparisons to the existing flux corrected transport (FCT) scheme in MPAS-Ocean.
In \firstRev{section 4 } the conclusions are discussed.
\firstRev{Finally the formulation of a local coordinate system tangent to the sphere in each element 
and the mapping between local and global coordinates are given in the appendix.}

\section{Formulation}

\subsection{Characteristic discontinuous Galerkin advection}

For a detailed description of the CDG scheme, the reader is referred to the previous article \cite{Lee16}.
Here we provide a brief overview of the formulation. The equation for the advection of a thickness-weighted 
tracer concentration is given as

\begin{equation}\label{eqn1.1}
\frac{\partial hq}{\partial t} + \nabla\cdot(\vec uhq) = 0
\end{equation}
where $h$ is the layer thickness, $q$ the tracer concentration and $\vec u$ the transport velocity.
\firstRev{Note that here the layer thickness is assumed to be constant with time. We will relax this assumption 
when discussing the implementation on the vertical ALE grid in section \ref{form_vert}. }
We begin by discretizing the domain into a set of $k$ contiguous elements, 
\firstRev{for which $h_k$ and $q_k$ are the discrete forms of the layer thickness and tracer concentration respectively. } 
Multiplying \firstRev{$h_kq_k$ }
by a set of $i$ test functions that vary in both space and time $\phi_{k,i}(\vec x,t)$, and expanding via 
the chain rule gives

\begin{equation}\label{eqn1.2}
\frac{\partial\phi_{k,i}h_kq_k}{\partial t} + \nabla\cdot(\phi_{k,i}\vec{u}h_{k}q_{k}) = 
\phi_{k,i}\Big(\frac{\partial h_kq_k}{\partial t} + \nabla\cdot(\vec{u}h_kq_k)\Big) + 
h_kq_k\Big(\frac{\partial\phi_{k,i}}{\partial t} + \vec{u}\cdot\nabla\phi_{k,i}\Big).
\end{equation}
\firstRev{Note that this formulation differs from the standard Galerkin formulation in that the tracer concentration
and not the full equation has been multiplied by the test function } \cite{Lee16}.
The first term on the right-hand side of \eqref{eqn1.2} is the discrete form of \eqref{eqn1.1} weighted by the test function
and so vanishes, and the second term represents the material derivative of the test functions $D\phi_{k,i}/Dt$. Integrating 
by parts over the element area $\Omega_k$ and between time levels $n$ and $n+1$, and then applying Gauss' theorem, the weak 
form is given as

\begin{multline}\label{eqn1.3}
\int_{\Omega_k}(\phi_{k,i}h_kq_k)^{n+1} - (\phi_{k,i}h_kq_k)^n\mathrm{d}\vec x + 
\int_{t^n}^{t^{n+1}}\int_{\partial\Omega_k}\phi_{k,i}\vec{u}h_{k'}q_{k'}\cdot\mathrm{d}\vec{s}\mathrm{d}t = \\
\int_{t^n}^{t^{n+1}}\int_{\Omega_k}h_kq_k\frac{D\phi_{k,i}}{Dt}\mathrm{d}\vec x\mathrm{d}t,
\end{multline}
\firstRev{where $k'$ denotes a set of elements within a local neighbourhood of element $k$
from which the flux is to be determined. }
The right-hand side may be taken as zero if the values of the test functions are constant along characteristics,
a condition also enforced in the ELLAM \cite{CRHE90,RC02} and semi-Lagrangian discontinuous Galerkin
\cite{GNQ14} schemes, such that

\begin{equation}\label{eqn1.4}
\frac{D\phi_{k,i}}{Dt} = 0.
\end{equation}
Equation (\ref{eqn1.4}) is satisfied via the introduction of a test function $\beta$ which varies with respect
to a Lagrangian coordinate in space and time $\vec\Gamma$ as

\begin{equation}\label{eqn1.5}
\phi_{k,i}(\vec x,t) = \beta_{k,i}(\vec\Gamma(\vec\xi(s),s))
\end{equation}
where $\vec\Gamma(\vec\xi(s),s)$ is constant with respect to the parametric variable $s$ along the
characteristic trajectory

\begin{equation}\label{eqn1.6}
\frac{\mathrm{d}\vec\xi}{\mathrm{d}s} = \vec u(\vec\xi(s),s)\qquad\vec\xi(t) = \vec x
\end{equation}
and $t$ is the point on $s$ where the boundary condition is applied, such that

\begin{equation}\label{eqn1.7}
\frac{\mathrm{d}\vec\Gamma(\vec\xi(s),s)}{\mathrm{d}s} = 0\qquad
\vec\Gamma(\vec x,t^{n+1}) = \vec\xi(t^{n+1}).
\end{equation}
Note that the boundary condition for \eqref{eqn1.7} follows from that for
\eqref{eqn1.6} such that for any $s = t$, $\Gamma(\vec\xi(t),t) = \Gamma(\vec x,t)$, with $t^{n+1}$ being the
specific time at which the boundary condition is applied.
Integrating (\ref{eqn1.6}) with respect to $s$ between $t$ and $t^{n+1}$
and recalling the boundary conditions on $\vec\xi(t)$ and $\vec\xi(t^{n+1})$ gives

\begin{equation}\label{eqn1.9}
\vec{\Gamma}(\vec x,t) = \vec{x} + \int_{t}^{t^{n+1}}\vec{u}(\vec\xi(s),s)\mathrm{d}s
\end{equation}
such that for any $(\vec{x},t)$ \eqref{eqn1.9} preserves the constant value of $\beta$ along
characteristics and hence the constant value of $\phi$ along those same characteristics.

Equation \eqref{eqn1.9} implies that the test functions \emph{arrive} at their static Eulerian coordinates
$\vec x$ at time level $n+1$ such that $\phi_{k,i}(\vec x,t^{n+1}) = \beta_{k,i}(\vec x)$,
with $\vec x$ being the location of $\vec\Gamma(\vec x,t^{n+1})$ as given in \eqref{eqn1.9}.
The values of the test functions at the same coordinate at the previous time level $n$ are then given as 
$\phi_{k,i}(\vec x,t^n) = \beta_{k,i}(\vec\Gamma(\vec x,t^n))$. 
As a practical matter, this means that if we wish to evaluate a test function $\phi_{k,i}(\vec x,t^n)$ at a given coordinate
$\vec x$ at a previous time level $n$ subject to \eqref{eqn1.4} then we may equivalently evaluate
the static test function $\beta_{k,i}(\vec\Gamma(\vec x,t^n))$ at its previous location by integrating 
\emph{forwards} with the Eulerian velocity field.

Unlike the test functions, $\phi_{k,i}(\vec x,t)$, there 
is no requirement that the trial functions be conserved along characteristics, and so may remain static, defined 
at the arrival locations of the characteristics at time level $n+1$ as given in \eqref{eqn1.9}. 
Ensuring that the mass matrix by which the solution coefficients are multiplied remained static in this fashion
motivated our choice of boundary conditions in \eqref{eqn1.6} and \eqref{eqn1.7}.
\firstRev{We represent the discrete form of the tracer concentration via an expansion of yet to be defined trial 
functions as }
\begin{equation}
q_{k}(\vec x,t) = \sum_jc_{k,j}(t)\beta_{k,j}(\vec x).
\end{equation}
This gives rise to the following linear system for the solution of the trial function coefficients $c_{k,j}^{n+1}$ 
in each element $k$ at the new time level $n+1$

\begin{equation}\label{eqn1.10}
\sum_j
\int_{\Omega_k}h_k\beta_{k,i}\beta_{k,j}\mathrm{d}\vec x c_{k,j}^{n+1} = \int_{\Omega_k}(\phi_{k,i}h_kq_k)^n\mathrm{d}\vec x - 
\int_{t^n}^{t^{n+1}}\int_{\partial\Omega_k}\phi_{k,i}\vec{u}h_{k'}q_{k'}\cdot\mathrm{d}\vec{s}\mathrm{d}t.
\end{equation}

Like the standard discontinuous Galerkin formulation, only the boundary fluxes are required to 
determine the solution of the tracer coefficients at the new time level, so no global mass matrix is required.
However, unlike the standard discontinuous Galerkin formulation, where the fluxes are determined 
via some Eulerian process, such as some form of averaging, upwinding or a Riemann solver, under the 
CDG formulation the edge fluxes must be evaluated by taking the area made by the edge as it is swept backward in time 
to its static Eulerian location from the previous time level and integrating the tracer mass over this area. This may 
be expressed as

\begin{equation}\label{eqn1.11}
\sum_j
\int_{\Omega_k}h_k\beta_{k,i}\beta_{k,j}\mathrm{d}\vec x c_{k,j}^{n+1} = \int_{\Omega_k}(\phi_{k,i}h_kq_k)^n\mathrm{d}\vec x - \\
\sum_e\sum_{k'}\int_{\Omega_{k,k',\Delta t}}(\phi_{k,i}h_{k'}q_{k'})^n\mathrm{d}\vec x
\end{equation}
where $\Omega_{k,k',\Delta t}$ is the intersection of the \emph{swept region} of the edge $e$ of element $k$
over time step $\Delta t$ and element $k'$.
The solution of (\ref{eqn1.11}) subject to (\ref{eqn1.4}) represents the characteristic discontinuous Galerkin formulation
for updating the tracer trial function coefficients $c_{k,j}^{n+1}$ at the new time level $n+1$. 

\begin{figure}[!hbtp]
\centering
\includegraphics[width=0.48\textwidth,height=0.40\textwidth]{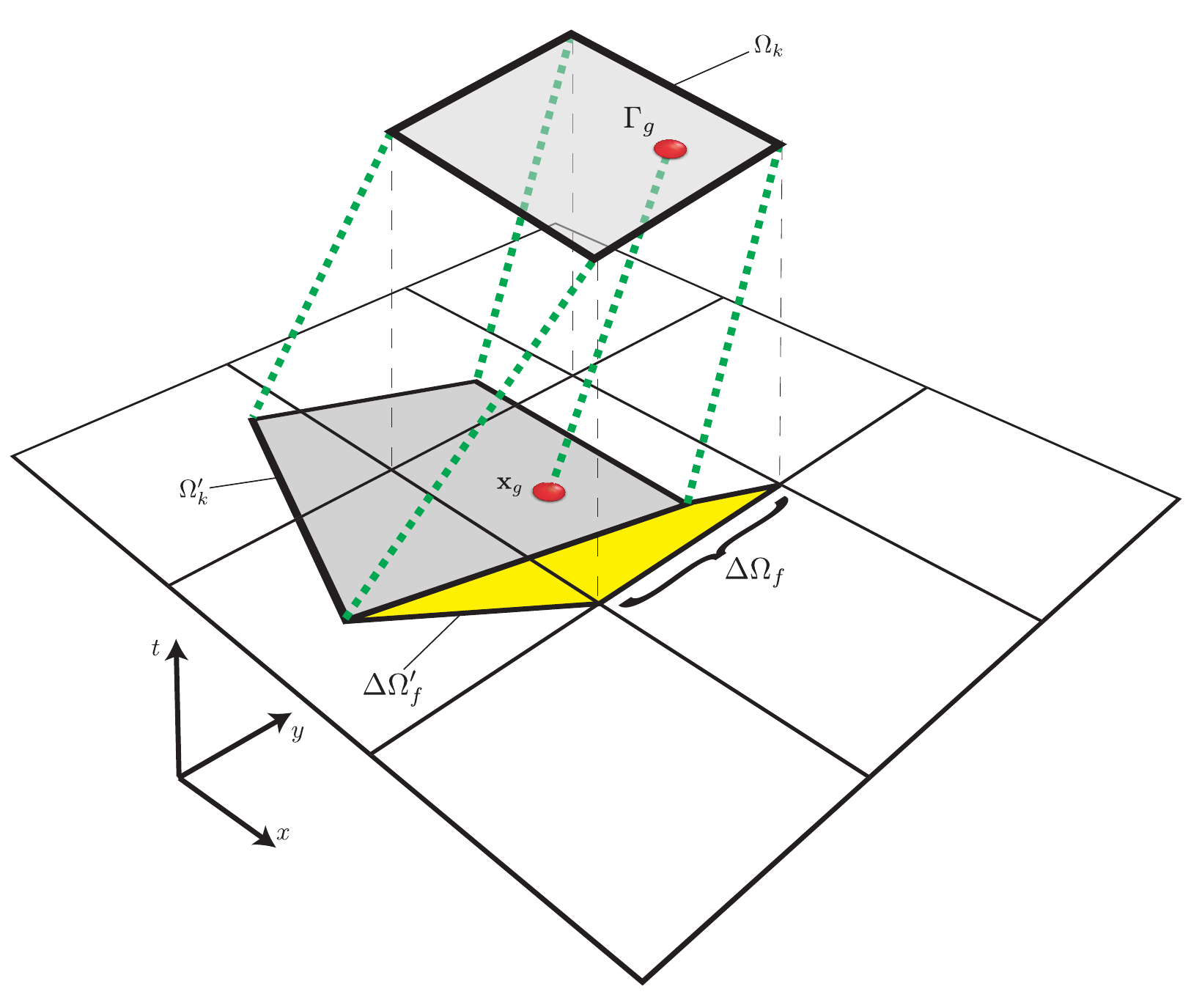}
\caption{Schematic of the flux computation for the CDG scheme. Edge vertices are integrated back in time to their
departure points in order to determine the swept region for an edge over a given time step. The volume of the tracer
over the swept region is then integrated with the quadrature points integrated forwards to arrival points where the 
test functions are evaluated in order to preserve the value of the test functions along characteristics.}
\label{Fig1}
\end{figure}

\subsection{Vertical CDG advection in ALE coordinates}\label{form_vert}

Before presenting the full 3D formulation of the CDG advection scheme we discuss the 1D vertical advection in
some detail. This is in order to present the subtleties required to conservatively apply the vertical advection 
scheme in the context of the moving ALE vertical grid \cite{Petersen15}.

The vertical advection of a tracer concentration $q(z,t)$ in a layer of varying thickness $h(z,t)$ is given as

\begin{equation}\label{eqn2.1}
\frac{\partial hq}{\partial t} + \frac{\partial(w-w_r)hq}{\partial z} = 0
\end{equation}
where $w$ is the vertical velocity and $w_r$ the velocity of the layer interface such that $w-w_r$ is the effective
transport velocity across the layer interface. 
\firstRev{As for the horizontal formulation in the preceeding section, we begin by multiplying the discrete form
of the tracer concentration $q_l$ by the $i$ test functions for level $l$, }
$\phi_{l,i}(z,t)$ 
\firstRev{and expanding via the chain rule as }

\begin{multline}\label{eqn2.2}
\frac{\partial\phi_{l,i}h_lq_l}{\partial t} + \frac{\partial\phi_{l,i}(w-w_r)h_{l}q_{l}}{\partial z} = 
\phi_{l,i}\Big(\frac{\partial h_lq_l}{\partial t} + \frac{\partial (w-w_r)h_lq_l}{\partial z}\Big) \\
+ h_lq_l\Big(\frac{\partial\phi_{l,i}}{\partial t} + (w-w_r)\frac{\partial\phi_{l,i}}{\partial z}\Big).
\end{multline}
\firstRev{Again, the first term on the right hand side vanishes as it is the 
discrete form of \eqref{eqn2.1} weighted by the test function. }
For the right-hand side to be zero the test function must move with velocity $w-w_r$ such that

\begin{equation}\label{eqn2.3}
\frac{\partial\phi_{l,i}}{\partial t} + (w-w_r)\frac{\partial\phi_{l,i}}{\partial z} = 0\qquad \phi_{l,i}(z,t^{n+1}) = \beta_{l,i}(z).
\end{equation}

The horizontal formulation of the CDG scheme given in equation \eqref{eqn1.10} assumes a constant layer thickness
$h$. For the rest of this article we relax this assumption in order to account for the varying layer thickness of 
the vertical ALE grid. Instead of solving for a thickness-weighted tracer concentration $hq(\vec x,t)$ subject to 
constant thickness, we must therefore integrate over the layer thickness in the vertical. Proceeding from this
representation we assume a trial function expansion as

\begin{equation}\label{eqn2.3.1}
q_l(z,t) = \sum_ja_{l,j}(t)\beta_{l,j}(z),
\end{equation}
and integrating with respect to space and time gives

\begin{multline}\label{eqn2.4}
\sum_j\int_{h_l^{n+1}}\beta_{l,i}\beta_{l,j}\mathrm{d}za_{l,j}^{n+1} = 
\int_{h_l^n}(\phi_{l,i}q_l)^n\mathrm{d}z - \\
\int_{\partial h_{l}^n}\int_{t^n}^{t^{n+1}}(w-w_r)\phi_{l,i}q_{l'}|_{l-} - (w-w_r)\phi_{l,i}q_{l'}|_{l+}\mathrm{d}t\mathrm{d}z
\end{multline}
where $l\pm$ denote the bottom and top of the layer interface respectively 
(with the vertical coordinate decreasing with layer index). 
Note that the vertical domains differ between time levels $n$ and $n+1$ as $h_l(z,t)$ evolves. 

\firstRev{We will assume }
that the departure regions for the top and bottom of layer \firstRev{$l$ } are determined at time level $n$,
\firstRev{noting that this is a particular choice for the representation of the flux terms and is not 
specifically required of the algorithm. }
We can express
the right hand side boundary terms as a sum over intersections between the swept region of the bottom and top 
interfaces over a time step $\Delta t$ of $l$ and the intersecting levels $l'$ (at time level $n$), 
$h_{l,l',\Delta t}^n$ as

\begin{equation}\label{eqn2.5}
\sum_j\int_{h_l^{n+1}}\beta_{l,i}\beta_{l,j}\mathrm{d}za_{l,j}^{n+1} = 
\int_{h_l^n}(\phi_{l,i}q_l)^n\mathrm{d}z - 
\int_{h_{l,l',\Delta t}^n}(\phi_{l,i}q_{l'})^n|_{l-} - (\phi_{l,i}q_{l'})^n|_{l+}\mathrm{d}z.
\end{equation}
\firstRev{Equation \eqref{eqn2.5} is the vertical analogue of the horizontal CDG advection equation \eqref{eqn1.11}. }
\secondRev{We note that the CDG scheme for vertical transport is fundamentally different from the widely used 
Lagrangian remap scheme \cite{Lin04}, in that the higher order moments are determined via a Galerkin projection 
of swept region fluxes onto these terms in the same fashion as the mean component, rather than being 
reconstructed from the mean values in the neighbouring layers. }

\firstRev{In \cite{Lee16} the horizontal CDG scheme was shown to be locally conservative. The vertical
scheme is also locally conservative, since the tracer \emph{mass} is the integral of the tracer concentration
$q_l$ over the element volume, $h_l^n$ for the 1D vertical scheme. Assuming a set of coefficients $b_i$ for 
the test functions, such that $\sum_ib_i\beta_i(z) = 1$, \eqref{eqn2.5} gives
\begin{multline}\label{eqn2.6}
\sum_ib_i\sum_j\int_{h_l^{n+1}}\beta_{l,i}\beta_{l,j}\mathrm{d}za_{l,j}^{n+1} = 
\sum_ib_i\int_{h_l^n}(\phi_{l,i}q_l)^n\mathrm{d}z \\
-
\sum_ib_i\int_{h_{l,l',\Delta t}^n}(\phi_{l,i}q_{l'})^n|_{l-} - (\phi_{l,i}q_{l'})^n|_{l+}\mathrm{d}z.
\end{multline}
Using the vertical analogue of \eqref{eqn1.5}, this simplifies to 
\begin{equation}\label{eqn2.7}
\int_{h_l^{n+1}}q_l^{n+1}\mathrm{d}z -
\int_{h_l^n}q_l^n\mathrm{d}z = 
-\int_{h_{l,l',\Delta t}^n}q_{l'}^n|_{l-} - q_{l'}^n|_{l+}\mathrm{d}z.
\end{equation}
The flux terms cancel with those from the neighbouring elements, since these are equal and opposite.
Since $q_l^n$ is a \emph{concentration}, the tracer \emph{mass} is the integral of $q_l^n$ over element $l$,
which is conserved between time levels.
}

\subsection{CDG advection in 3D}

Having presented the formulation of the CDG advection scheme independently in the horizontal and the
vertical, we now proceed to the formulation of the full 3D scheme.
The 3D advection of the tracer concentration is given as

\begin{equation}\label{eqn3.1}
\frac{\partial hq}{\partial t} + \nabla\cdot(\vec uhq) + \frac{\partial(w-w_r)hq}{\partial z} = 0.
\end{equation}

We assume a modal Taylor series test function expansion in both the horizontal and vertical dimensions as

\begin{multline}\label{eqn3.2.1}
\beta_{k,l}(x,y,z) = \sum_i\beta_{k,l,i}(x,y,z) = 1 + \frac{1}{\Delta x}(x - \overline{x}) + \frac{1}{\Delta y}(y - \overline{y}) + \\
\frac{1}{2\Delta x^2}(x^2 - \overline{x^2}) + \frac{1}{\Delta x\Delta y}(xy - \overline{xy}) + \frac{1}{2\Delta y^2}(y^2 - \overline{y^2}) + ... + \\
\frac{1}{\Delta z}(z - \overline{z}) + \frac{1}{\Delta z^2}(z^2 - \overline{z^2}) + ...
\end{multline}
\firstRev{where $k$ and $l$ are the element indices in the horizontal and vertical dimensions respectively, and }
$\Delta x$, $\Delta y$ and $\Delta z$ are length scales of the element in $x$, $y$ and $z$ respectively 
by which the terms are normalized in order to keep them $\mathcal{O}(1)$. The terms denoted by the overbars are 
the mean components defined as

\begin{equation}\label{eqn1.9.3}
\overline{x^my^n} = \int_{\Omega_k}x^my^n\mathrm{d}\vec x\mathrm{d}z\qquad \overline{z^m} = \int_{\Omega_k}z^m\mathrm{d}\vec x\mathrm{d}z
\end{equation}
which are removed from the higher order terms to ensure that they remain massless, such that a slope limiter may be
applied to these without loss of conservation. We also assume a trial function expansion for $q$ in each element as

\begin{equation}\label{eqn3.2.2}
q_{k,l}(\vec x,z,t) = \sum_jc_{k,l,j}(t)\beta_{k,l,j}(\vec x,z).
\end{equation}

Recalling the horizontal and vertical formulations of the CDG scheme as given in \eqref{eqn1.11} and 
\eqref{eqn2.5} respectively gives

\begin{multline}\label{eqn3.3}
\sum_j\int_{h_{k,l}^{n+1}}\int_{\Omega_{k,l}}\beta_{k,l,i}\beta_{k,l,j}\mathrm{d}\vec x\mathrm{d}zc_{k,l,j}^{n+1} =
\int_{h_{k,l}^n}\int_{\Omega_{k,l}}(\phi_{k,l,i}q_{k,l})^n\mathrm{d}\vec x\mathrm{d}z - \\
\sum_e\frac{V_e^{con}}{V_e^{cdg}}
\sum_{k'}\int_{h_{k,k',l}^n}\int_{\Omega_{k,k',l,\Delta t}}(\phi_{k,l,i}q_{k',l})^n\mathrm{d}\vec x\mathrm{d}z - \\
\sum_{l'}\int_{h_{k,l,l',\Delta t}^n}\int_{\Omega_{k,l,l'}}(\phi_{k,l,i}q_{k,l'})^n|_{l-} - (\phi_{k,l,i}q_{k,l'})^n|_{l+}\mathrm{d}\vec x\mathrm{d}z,
\end{multline}
\firstRev{where $k'$ and $l'$ are the set of elements in the horizontal and vertical respectively which intersect 
with the swept region of element $k$, $l$ over time step $\Delta t$. }
Note that the layer thicknesses are determined separately from the continuity equation at integer time steps.
\secondRev{While the fluxes have been partitioned into their horizontal and vertical components in \eqref{eqn3.3},
these are still applied at the same time level. As such there is no time splitting of the horizontal and vertical
advection operators, and both the horizontal and vertical flux terms project onto the three dimensional solution 
of the tracer concentration. }

The horizontal fluxes have been normalized by the factor $V_e^{con}/V_e^{cdg}$. This represents the ratio
of the volume fluxed across edge $e$ as determined from the (single moment, finite volume) continuity equation 
and that swept across the same edge using the CDG algorithm. The volume fluxed across the edge from the 
continuity equation is centered in space, given as

\begin{equation}
V_e^{con} = 0.5(h_{e-} + h_{e+})\Delta t\vec u\cdot\vec nd_e
\end{equation}
where $d_e$ is the width of the edge and $h_{e-}$, $h_{e+}$ are the thicknesses of the elements either 
side of the edge. Unlike the continuity equation volume fluxes, the CDG swept region flux is upwinded, 
and may be given by a piecewise constant integration of the layer thickness over the swept regions as

\begin{equation}
V_e^{cdg} = \sum_{k'}\int_{\Omega_{k,k',\Delta t}}h_{k'}^n\mathrm{d}\vec x.
\end{equation}
This scaling of the edge fluxes ensures that the implicit volume fluxes of the CDG scheme are
\emph{consistent} with respect to the explicit volume fluxes from the continuity equation.
This procedure for enforcing consistency is much simpler than that required for remap schemes based
on a multi-moment representation of the continuity equation \cite{Lauritzen16}.
\secondRev{The relative simplicity of the consistency fix presented here is based on i. the fact that
the layer thickness is represented via a single moment in each element in keeping with the finite volume
formulation of the continuity equation and ii. the formulation of the CDG scheme in flux form rather
than remap form, such that the edges and not the elements are traced back along characteristics. }
It is worth noting that while the CDG 
scheme uses an upwinded flux, the volume flux used by the continuity equation is centered. This is necessary
since the continuity equation must account for both the left and right gravity wave solutions of the
linearized shallow water equations, whereas no such restriction is required for the CDG advection
of passive tracers.

While the vertical and horizontal fluxes are evaluated separately in order to avoid the need to evaluate 
swept region intersection in three dimensions, \eqref{eqn1.4} 
\firstRev{is still }
satisfied in 3D as

\begin{equation}
\frac{\partial\phi}{\partial t} + \vec u\cdot\nabla\phi + (w-w_r)\frac{\partial\phi}{\partial z} = 0\qquad
\phi(\vec x,z,t^{n+1}) = \beta(\vec x,z).
\end{equation}
\firstRev{This is to ensure that the effects of both the horizontal and vertical fluxes are accounted for
in the advection of the test functions.}

\section{Results}

The CDG scheme has previously been verified via convergence studies for analytic solutions on planar
quadrilateral and hexahedral grids \cite{Lee16}. Here we present comparisons to the existing FCT 
scheme in MPAS-Ocean \cite{Ringler13} for both passive advection on the sphere \cite{LSPT12} and 
a suite of idealized test cases with full ocean dynamics \cite{Ilicak12, Petersen15}.

While the CDG scheme may be run to arbitrarily high order accuracy (provided that a quadrature 
rule \secondRev{and a mapping } exist to integrate the curvature of the pre image of each edge to 
the desired order of accuracy), in each of the test cases presented here we use a linear basis in 
each dimension. 
\secondRev{The principle reason for this is that while higher-order mappings from the sphere
to elements in the plane are relatively straight forward to construct for methods on quadrilateral 
tensor product elements \cite{GNQ14,Lauritzen16}, a method for doing this on Voronoi elements with
an arbitrary number of sides is not known to the authors. }
The second order coefficients are slope limited using either a 3D implementation 
of the vertex based Barth-Jespersen limiter \cite{Kuzmin10} or simplified WENO method 
which uses the basis functions within each element to determine the smoothing coefficients \cite{ZS13,GNZ16}.

\secondRev{As reported previously \cite{Lee16}, the CDG scheme, like other semi-Lagrangian methods, is
unconditionally stable with time step. Provided that the halo size used for the parallel decomposition of 
the domain is sufficiently large and the characteristics for a given edge do not cross one another 
\cite{LR05}, then the scheme may be run with any CFL number. While the scheme has been run on the sphere 
with CFL numbers up to 2.5, in the results presented here we limit ourselves to time steps equal to those
used by the dynamics $(\mathrm{CFL} < 1)$. This is because it is difficult to robustly preserve consistency 
with respect to the continuity equation with larger CFL numbers, due to variations in layer thickness. The 
scheme requires a single halo update at the end of each time step in order to ensure that the tracer fields 
on the boundaries are consistent for each processor, however we note that this halo update must include all 
the moments and not just the mean components as is the case for the FCT scheme.}

\subsection{Test case 1: passive advection}

The first test case involves the passive advection of a tracer field 
\firstRev{with a Gaussian initial distribution } within a deformational shear flow
on the sphere \cite{LSPT12}. \firstRev{Note that this configuration serves only to test the 
convergence of errors for the horizontal scheme on the sphere. } The $L_2$ errors are computed 
after 12 days when the tracer field has
returned to its original position, for both the CDG scheme and the existing FCT scheme 
\cite{BB73,Zalesak79} within MPAS-Ocean.

\begin{figure}[!hbtp]
\centering
\includegraphics[width=0.51\textwidth,height=0.38\textwidth,valign=c]{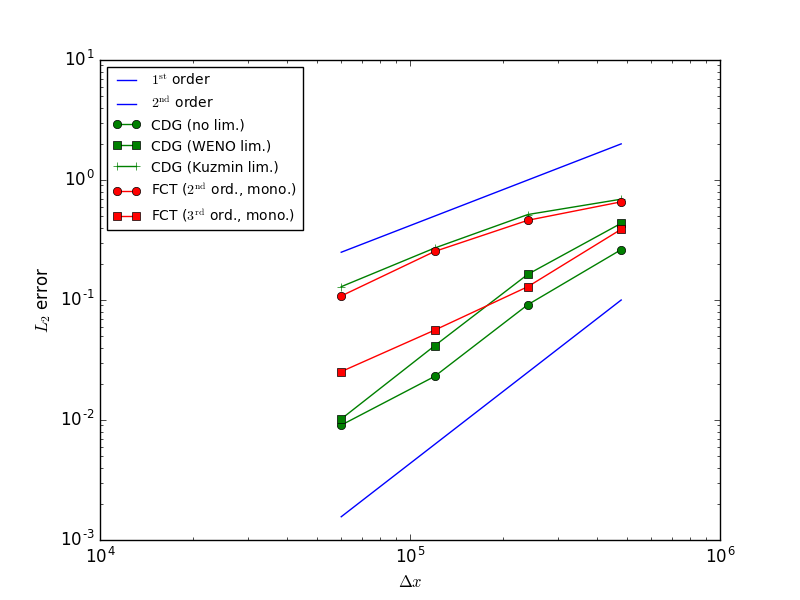}
\includegraphics[width=0.48\textwidth,height=0.28\textwidth,valign=c]{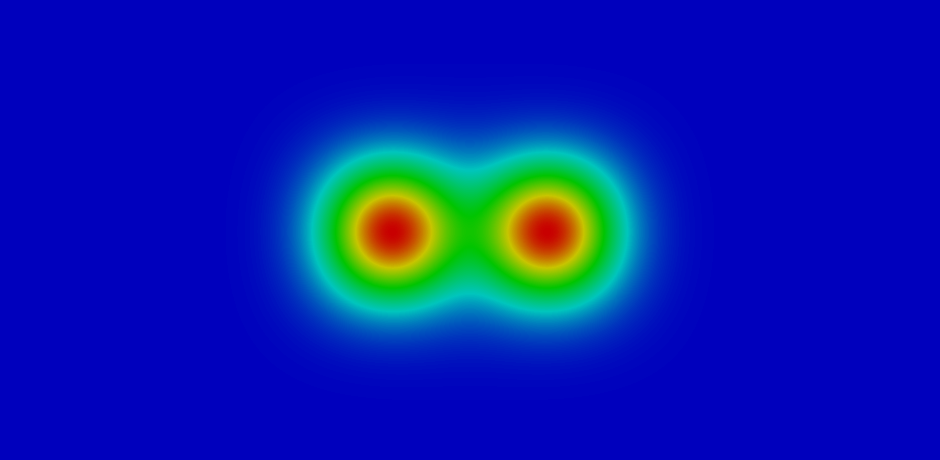}
\caption{Left: $L_2$ errors for the CDG and FCT transport schemes for passive advection for the 
deformational shear flow on the sphere test case. Right: Initial condition.}
\label{Fig2}
\end{figure}

\begin{figure}[!hbtp]
\centering
\includegraphics[width=0.48\textwidth,height=0.28\textwidth]{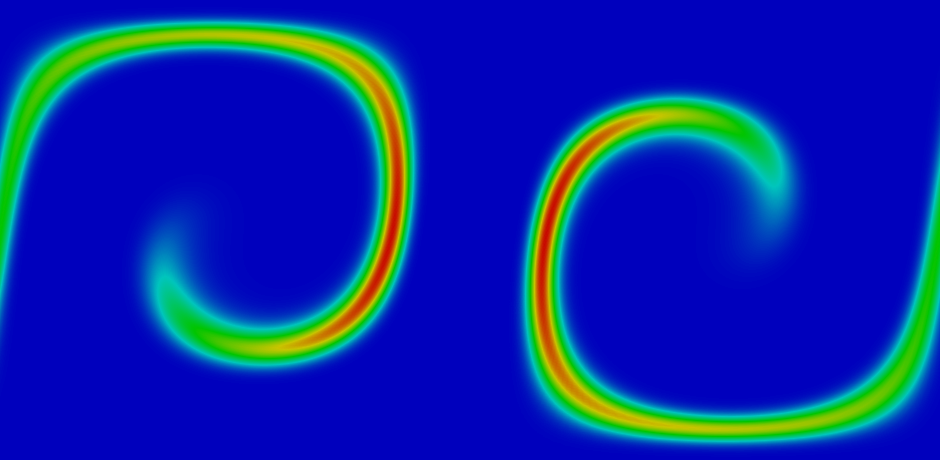}
\includegraphics[width=0.48\textwidth,height=0.28\textwidth]{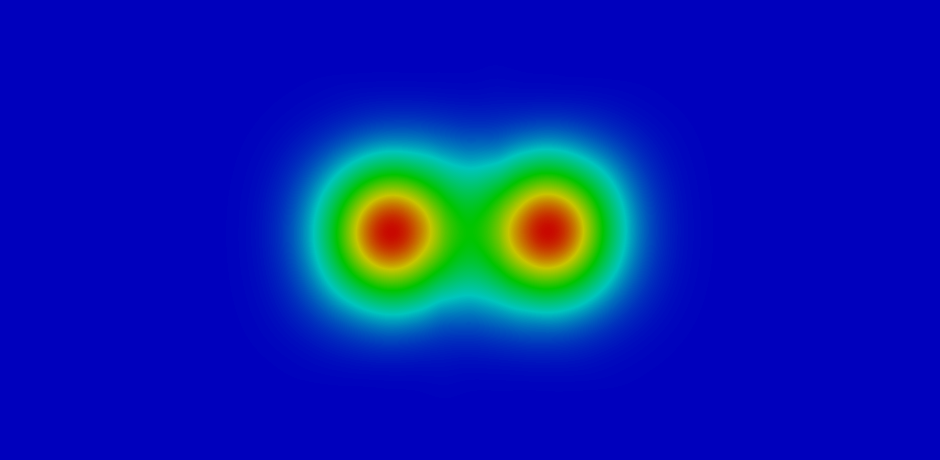}
\caption{Deformational shear flow advection on the sphere using CDG after 6 hours (left) and 
12 hours (right).}
\label{Fig3}
\end{figure}

As can be seen from fig. \ref{Fig2}, the CDG scheme for passive tracer transport on the sphere compares
favorably to the existing FCT scheme with MPAS ocean. The second order CDG scheme displays an error 
convergence rate superior to the third order FCT scheme in both the unlimited and WENO limited case.
For the subsequent active ocean test cases presented below we use the vertex based slope limiter 
\cite{Kuzmin10} as the WENO limiter is not able to preserve strict monotonicity.

\subsection{Test case 2: lock exchange}

\begin{figure}[!hbtp]
\centering
\includegraphics[width=0.51\textwidth,height=0.38\textwidth,valign=c]{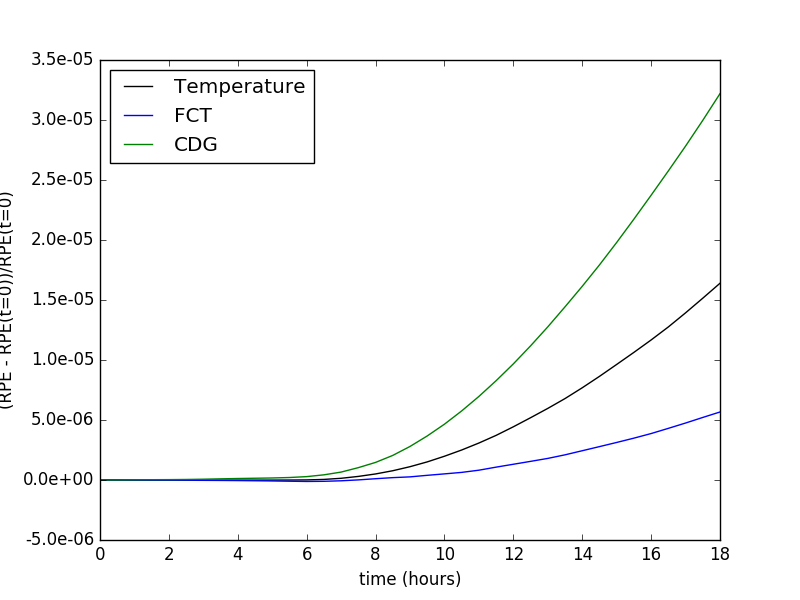}
\includegraphics[width=0.48\textwidth,height=0.36\textwidth,valign=c]{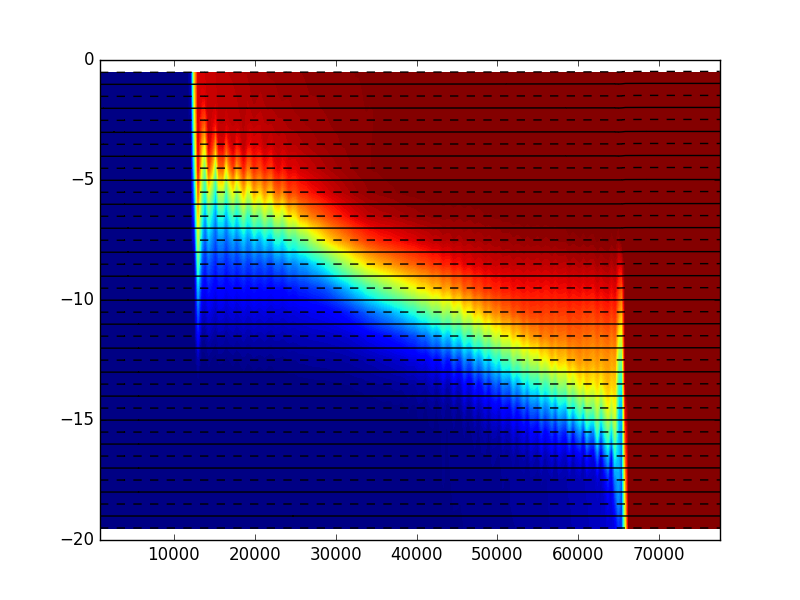}
\caption{Left: resting potential energy (RPE) with time for the lock exchange test case.
Right: temperature field (active tracer) at 18 hours. Values are between $5^{\circ}$ C and $30^{\circ}$ C.}
\label{Fig4}
\end{figure}

\begin{figure}[!hbtp]
\centering
\includegraphics[width=0.48\textwidth,height=0.36\textwidth]{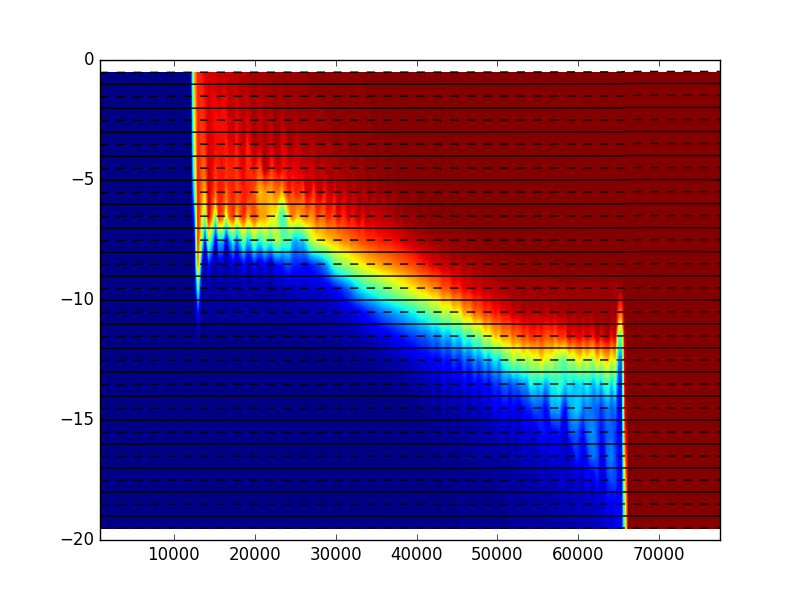}
\includegraphics[width=0.48\textwidth,height=0.36\textwidth]{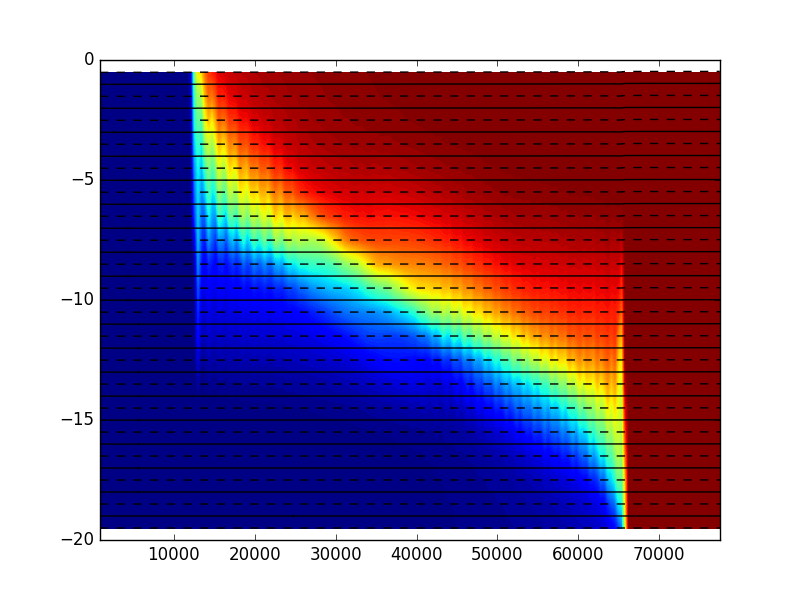}
\caption{Passive tracer after 18 hours using second order FCT (left) and CDG (right) advection.
Values are between $3^{\circ}$ C and $30^{\circ}$ C.}
\label{Fig5}
\end{figure}

For the second test case, an initial temperature distribution of $T=5^\circ$C on the left and $T=30^\circ$C
on the right side of a box generates a pressure gradient that drives a flow of sinking cool fluid along
the bottom to the right and rising warm fluid along the top to the left. The model uses a linear temperature ($T$)
dependent equation of state in order to determine the density $\rho $ of the form 
$\rho = 1000.0 - 0.2(T - 5.0)$, from which the pressure is derived via hydrostatic balance.
The initial temperature and passive tracer fields are given as
$q(y) = 5.0 + 12.5(1.0 + 2.0^{-4}\tanh(y - y_0))$. 
\firstRev{The model has periodic boundary conditions in the $x$ dimension with just 16 hexagonal elements along 
the periodic channel such that the dynamics are weak in this dimension. The ALE grid is configured such that the 
vertical height of the elements within each column stretch uniformly with perturbations in sea surface height, 
which are barely perceptible in figures \ref{Fig4} and \ref{Fig5}. }
Details of the specific geometry and model configuration can be found in \cite{Petersen15}.

The resting potential energy (RPE) is determined in order to quantify the amount of spurious
vertical mixing of the CDG scheme with respect to the existing FCT scheme for passive tracer
transport. The RPE is computed by reordering all the elements of the domain, by descending 
order of density $\rho$, into a single one dimensional column and then integrating this reorderd 
density $\rho^*$ down the column as $\mathrm{RPE} = \int_{\Omega}g\rho^*z\mathrm{d}V$
\cite{Ilicak12, Petersen15}. The RPE is computed both for the CDG and FCT schemes for passive
tracer transport, as well as for the FCT derived temperature as a reference. While the passive
FCT tracers are integrated using a first order forward Euler scheme, the active temperature
is solved using an iterated shooting method \cite{Ringler13} in order to derive a second order 
scheme which is necessary in order to ensure model stability.


As can be seen from fig. \ref{Fig4} the amount of spurious vertical mixing as measured by the RPE 
is significantly higher for the CDG scheme than either the passive or active FCT tracers. This
suggests that while the current limiting approach preserves monotonicity it is excessively diffusive
and new limiting strategies should be explored. This can also be seen by inspecting the tracer
fields at the final time, as given in figs. \ref{Fig4} and \ref{Fig5}.

\subsection{Test case 3: overflow}

The third test case also involves a horizontal step in temperature (between $10^\circ$C and $20^\circ$C),
however in this case a vertical step in topography is also included, such that as the cool fluid is driven 
down and rightward it also sinks down the topographic slope. The model configuration is similar to that for
the lock exchange test case, \firstRev{with periodic boundary conditions in the $x$ dimension and uniform 
stretching of the elements in the vertical in proportion to the perturbations
in sea surface height due to the barotropic mode, } and can be found in \cite{Petersen15}. This test serves 
to demonstrate that the CDG scheme remains stable for large vertical motions along steep topography.

\begin{figure}[!hbtp]
\centering
\includegraphics[width=0.48\textwidth,height=0.36\textwidth]{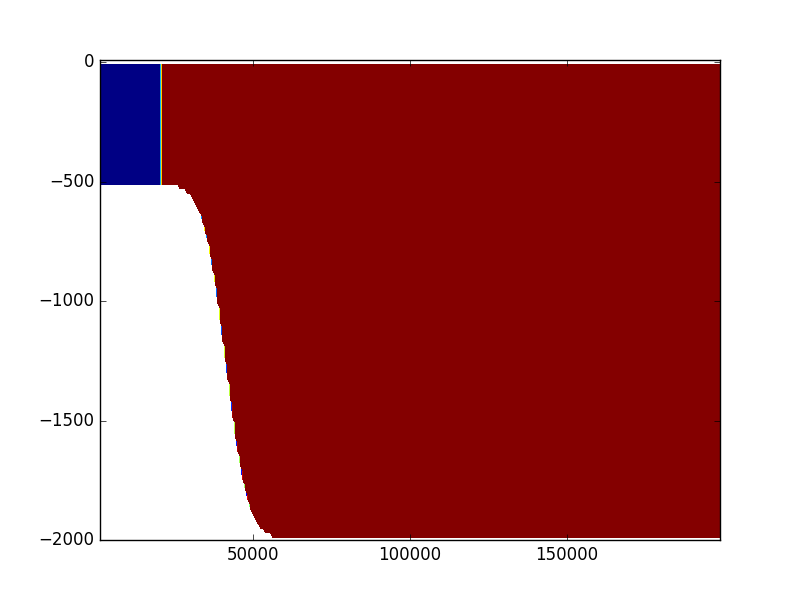}
\includegraphics[width=0.48\textwidth,height=0.36\textwidth]{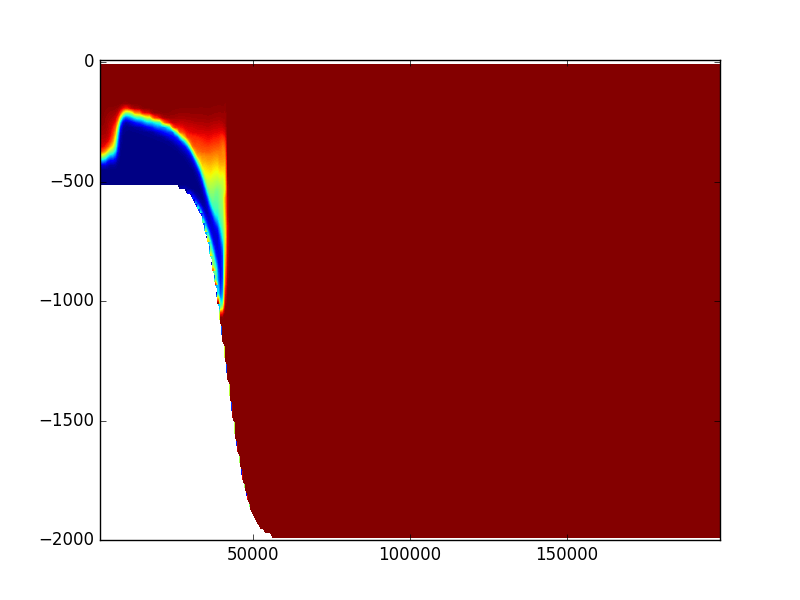}
\caption{Tracer field $q$ at $t = 0$ hours (left), $t = 5$ hours (right) for the 
overflow test case \cite{Petersen15}. Values are between $10^{\circ}$ C and $20^{\circ}$ C.}
\label{Fig6}
\end{figure}
\begin{figure}[!hbtp]
\centering
\includegraphics[width=0.48\textwidth,height=0.36\textwidth]{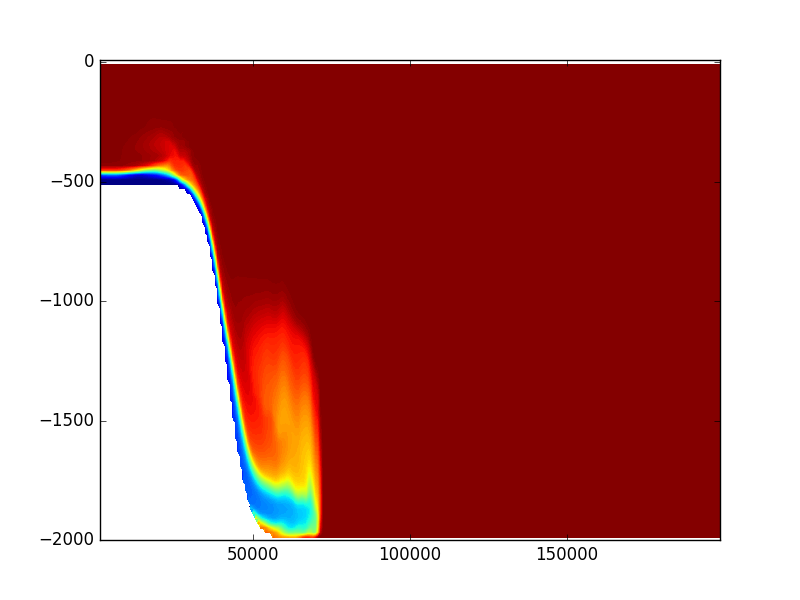}
\includegraphics[width=0.48\textwidth,height=0.36\textwidth]{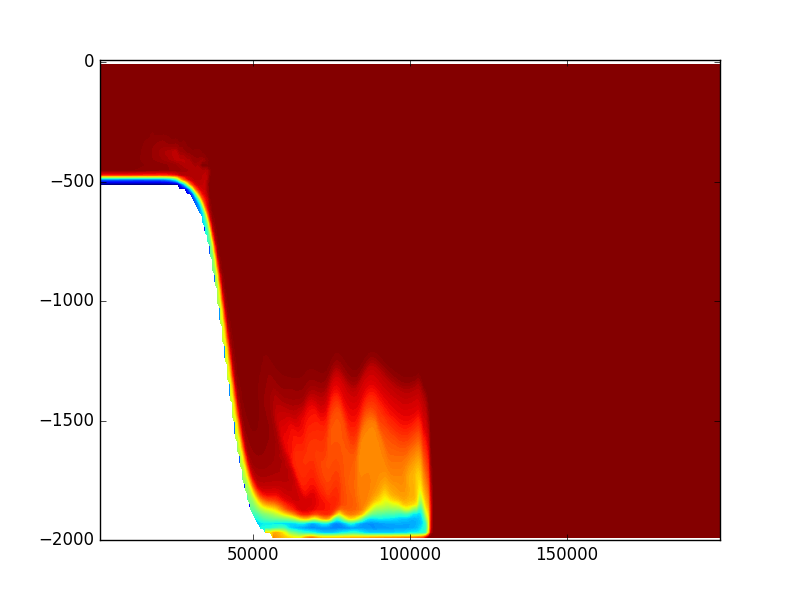}
\caption{Tracer field $q$ at $t = 10$ hours (left), $t = 15$ hours (right) for the 
overflow test case \cite{Petersen15}. Values are between $10^{\circ}$ C and $20^{\circ}$ C.}
\label{Fig7}
\end{figure}

Figures \ref{Fig6} and \ref{Fig7} show that the CDG scheme appropriately represents the flow of 
passive tracers along a steeply varying slope. The results are qualitatively similar to those 
previously published for the FCT scheme \cite{Petersen15}. However the quality of the results here
are somewhat misleading since the overly diffusive solution as demonstrated for the previous test
case serve to stabilize the strong downward motion, whereas the less diffusive FCT scheme requires
an explicit vertical viscosity in order to suppress numerical instabilities.

\subsection{Test case 4: baroclinic channel}

The baroclinic channel test case is initialized with vertically sloping isotherms in the meridional direction, 
as well as a sinusoidal temperature profile in the plane which drive the formation of baroclinically unstable
eddies \cite{Petersen15}. Unlike the previous test cases, the baroclinic channel configuration allows for 
significant motion in all dimensions, and so more fully supports the evolution of nonlinear momentum transport. 
Properly resolving  the formation and transport of the resultant eddies presents a significant challenge for 
the CDG scheme due to the additional numerical diffusion introduced by the slope limiting.

\begin{figure}[!hbtp]
\centering
\includegraphics[width=0.25\textwidth,height=0.75\textwidth]{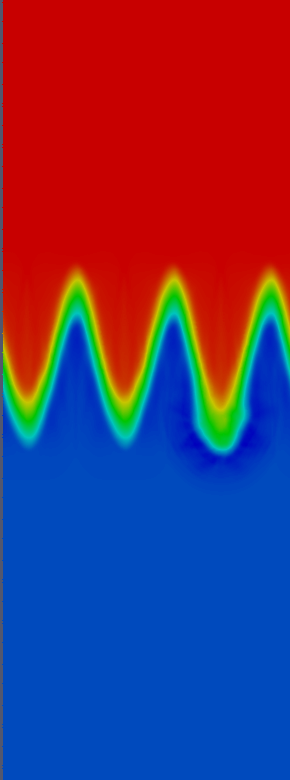}
\includegraphics[width=0.25\textwidth,height=0.75\textwidth]{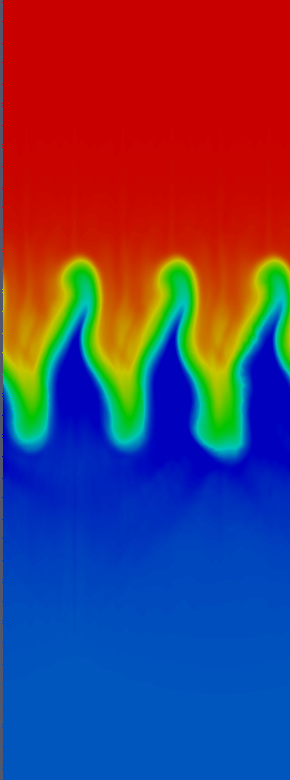}
\includegraphics[width=0.25\textwidth,height=0.75\textwidth]{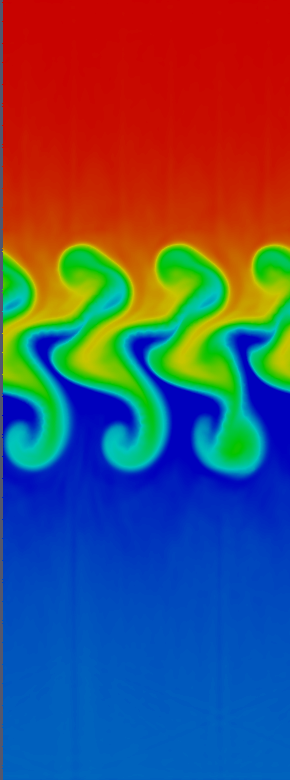}
\caption{Tracer field $q$ at $t = 0$ days (left), $t = 6$ days (center),
and $t = 12$ days (right) for the baroclinic channel test case \cite{Petersen15}
with channel length $400km$ and resolution $\Delta x = 2.5km$.}
\label{Fig8}
\end{figure}

Given a linear surface to bottom temperature profile with a difference of $3^\circ$C, a depth of $H=1000 m$, a reference
density of $\rho_0 = 1000 km/m^3$ and a Coriolis parameter of $f = -1.2\times 10^{-4}s^{-1}$, and recalling
the linear equation of state, the first baroclinic deformation radius is given as $L_d = \sqrt{g/\rho_0\partial p/\partial z}H/\pi f \approx 6.4 km$.
With a horizontal resolution of $\Delta x = 1 km$ the resultant eddies are just within the range of resolved scales 
of the simulation. Therefore the maintained presence of these eddies in the tracer field is highly sensitive to any
spurious numerical diffusion due to excessive slope limiting. For this reason the eddy field is observed to be significantly
muted in the tracer field compared to the FCT advected temperature field \cite{Petersen15}, as shown in fig. \ref{Fig8}.

\subsection{Test case 5: global ocean}

The final test case involves the spin up of a global ocean domain with a resolution of $120 km$.
It is initialized with climatological temperature and salinity and is driven by a temporally constant
surface wind profile. In order to successfully implement the CDG scheme in this configuration the departure and 
quadrature point integrations and the swept region intersections are performed on the sphere assuming great circle 
arcs for each of the edges, and the updated tracer coefficients for each element are solved by projecting each cell
into the plane as described in the appendix.

\begin{figure}[!hbtp]
\centering
\includegraphics[width=0.80\textwidth,height=0.40\textwidth]{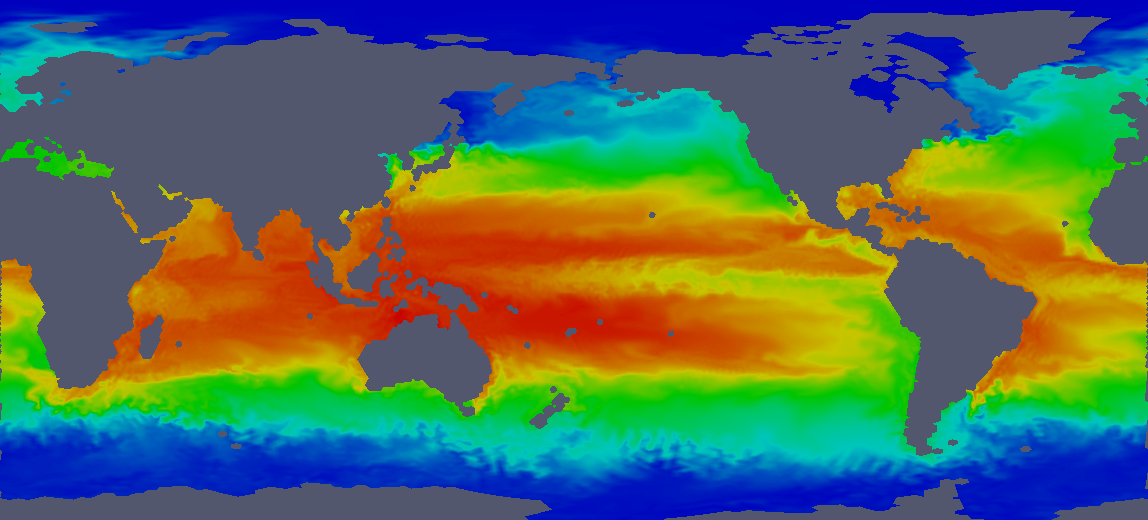}
\includegraphics[width=0.80\textwidth,height=0.40\textwidth]{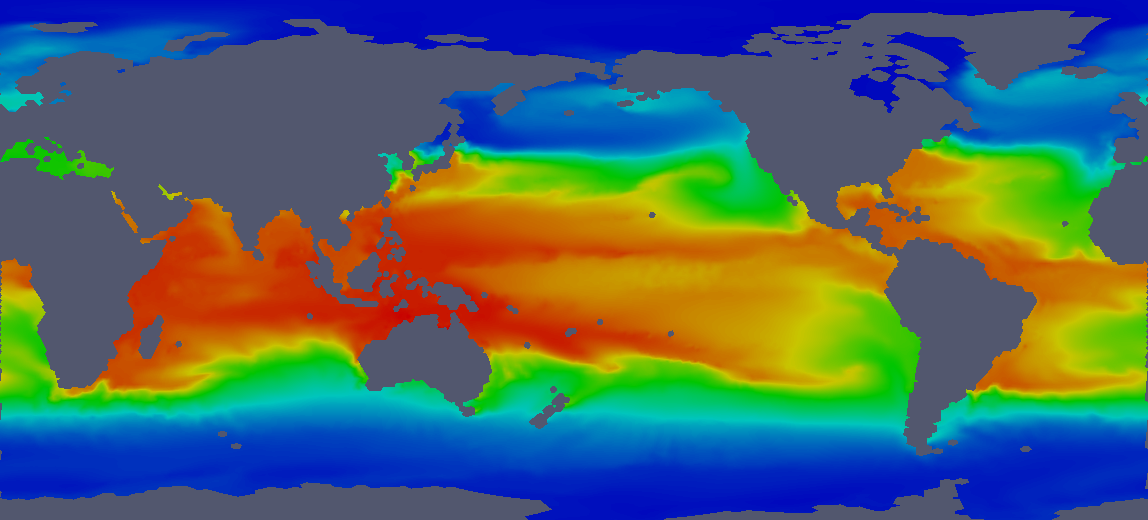}
\caption{Passive tracer  (FCT, $2^{\mathrm{rd}}$ order, top) and passive tracer
(CDG, $2^{\mathrm{nd}}$ order, bottom) for global ocean after 9 months. Color bars range
from $-1.8^{\circ}$ C to $30^{\circ}$ C.}
\label{Fig9}
\end{figure}

As can be seen in fig. \ref{Fig9} the CDG solution for the temperature field looks
broadly similar to the FCT solution. However the equatorial temperature is weaker and the 
secondary features such as western boundary currents are less pronounced, as seen in the Gulf Stream
and Kuroshio current.

\section{Conclusion}

We have presented an implementation of the characteristic discontinuous Galerkin (CDG) tracer transport scheme 
within the MPAS-Ocean model for both horizontal advection on an unstructured Voronoi grid and vertical advection 
on a temporally varying arbitrary Lagrangian-Eulerian (ALE) grid. The scheme has been used to model passive 
advection for a suite of idealized test cases, using the resting potential energy (RPE) as a measure of spurious
vertical mixing. Consistency between the implicit volume flux of the CDG scheme and the explicit volume
flux from the dynamics is enforced via a renormalization of the edge fluxes with respect to the volume
fluxes from the continuity equation.
Since the layer thickness is piecewise constant, this process is much simpler than the consistency fixers 
required for higher order representations of the thickness \cite{Lauritzen16}.

While the results compare favorably to the FCT scheme for passive advection with prescribed velocity on the 
surface of the sphere, the absence of a better limiting scheme and the need to preserve strict monotonicity
leads to excessive \firstRev{diffusion}, which significantly degrades the CDG solution in 3D. 
\secondRev{A significant issue with the slope limiting approach used here is that the same limiting
coefficient is used for moments in all dimensions for a given element, irrespective of which dimension has 
a non-monotone solution. }
As future work an improved limiter is required, perhaps using the recently developed anisotropic approach 
\cite{AKKR17},
\secondRev{so as to ensure that the limiting coefficients used to preserve monotone solutions in the vertical
are not projected onto the horizontal moments and vice versa. }

\section*{Acknowledgements}

The authors are grateful to Drs Mark Taylor, Peter Bosler and Andrew Bradley at Sandia National 
Laboratory for many enlightening discussions concerning the formulation of the CDG algorithm. We
would particularly like to thank Dr. Bradley for supplying the library for computing polygon 
intersections on the sphere.

We also thank Dr. Bill Lipscombe for supplying the method for performing tangent projections from the
sphere, the two anonymous reviewers, who's helpful comments greatly improved the clarity of this article.

We would also like to acknowledge the support of LANL 
Institutional Computing. This work was carried out under the auspices of the National Nuclear Security 
Administration of the U.S. Department of Energy at Los Alamos National Laboratory under Contract No. 
DE-AC52-06NA25396. This research was supported by the Office of Science (BER), U.S. Department of 
Energy. Los Alamos Report LA-UR-17-22608.

%\section*{References}

%\bibliographystyle{elsarticle-num} 
%\bibliography{references}

%\begin{thebibliography}{99}

%\end{thebibliography}

\section*{Appendix: Implementation on the sphere}

The implementation of the CDG algorithm on the sphere presents significant challenges as the trial 
and test functions are defined with respect to a local planar coordinate system. Here we present
a length preserving transformation between great circle arcs on the sphere and a plane tangent to
the element centers. The mapping between the sphere and the tangent plane is similar to that used 
in the MPAS-Seaice model \cite{Turner18}, with the distinction being that here lengths are preserved 
between the sphere and the tangent plane. 

The coordinates on the surface of a sphere of radius $R$ may be given in Cartesian space as

\begin{equation}\label{eqn4.1}
\vec x = [x\vec e_x, y\vec e_y, z\vec e_z].
\end{equation}
In order to construct an orthogonal local coordinate system for a tangent plane

\begin{equation}\label{eqn4.2}
\vec a = [a\vec e_a, b\vec e_b, c\vec e_c]
\end{equation}
we begin by defining the normalized vertical coordinate of the tangent plane as the being parallel
to the radial vector from the center of the sphere to the origin of the tangent plane $[x_o, y_o, z_o]^T$ such that
with respect to the global coordinates $\vec x$

\begin{equation}\label{eqn4.3}
\vec e_c = \Big[\frac{x_o}{R}, \frac{y_o}{R}, \frac{z_o}{R}\Big]
\end{equation}
where $R = \sqrt{x_o^2 + y_o^2 + z_o^2}$ is the constant radius of the sphere.
A unit vector denoting the $\mathbf{a}$ direction of the tangent plane is then given by projecting
$\vec e_c$ into the $x-y$ plane by setting the $z$ component to $0$ and rotating $\pi/2$ radians anti clockwise as

\begin{equation}\label{4.4}
\vec e_a = \Big[\frac{-y_o}{r}, \frac{x_o}{r}, 0\Big]
\end{equation}
where $r = \sqrt{x_o^2 + y_o^2}$ is the minor radius of the sphere at $z_o$.
The unit vector denoting the $\vec b$ direction of the tangent plane is then given by a right hand rule cross
product of these as

\begin{equation}\label{4.5}
\vec e_b = \vec e_c\times\vec e_a = \Big[\frac{-x_oz_o}{rR}, \frac{-y_oz_o}{rR}, \frac{r}{R}\Big].
\end{equation}
The transformation matrix from the local to the global coordinate frame is given as

\begin{equation}\label{4.6}
A = \begin{bmatrix}\vec e_a \\\vec e_b \\\vec e_c\end{bmatrix} =
\begin{bmatrix}
\frac{-y_o}{r}     & \frac{x_o}{r}      & 0              \\
\frac{-x_oz_o}{rR} & \frac{-y_oz_o}{rR} & \frac{r}{R}    \\
\frac{x_o}{R}      & \frac{y_o}{R}      & \frac{z_o}{R}  \\
\end{bmatrix}.
\end{equation}
With the mapping given as

\begin{equation}\label{4.7}
\vec a = A(\vec x - \vec x_o).
\end{equation}
Since $A$ represents a linear rotation of coordinate systems via unit vectors, the inverse mapping is simply by the
transpose $A^{-1} = A^T$ such that

\begin{equation}\label{eqn4.8}
\vec x = A^T\vec a + \vec x_o.
\end{equation}
A coordinate in the tangent plane with respect to the $\mathbf{a}$ coordinate system is given as $P_{\vec a} = [a,b,0]^T$,
where $c = 0$ since $P_{\vec a}$ is on the plane. This coordinate is transformed into global coordinates via (\ref{eqn4.8}) as
\begin{equation}
P_{\vec x} = 
\begin{bmatrix}
\frac{-y_oa}{r} - \frac{x_oz_ob}{rR} + x_o\\
\frac{x_oa}{r} - \frac{y_oz_ob}{rR} + y_o\\
\frac{rb}{R} + z_o\\
\end{bmatrix}.
\end{equation}
If $P_{\vec x}$ is to be projected down onto the sphere in a direction normal to the tangent plane then is must be shifted
a distance $\sqrt{R^2 - a^2 - b^2} - R$ in the $\mathbf{e}_c$ direction, such that this point is given as

\begin{equation}\label{eqn4.9}
Q_{\vec x} = 
\begin{bmatrix}
\frac{-y_oa}{r} - \frac{x_oz_ob}{rR} + \frac{\sqrt{R^2 - a^2 - b^2}x_o}{R}\\
\frac{x_oa}{r} - \frac{y_oz_ob}{rR} + \frac{\sqrt{R^2 - a^2 - b^2}y_o}{R}\\
\frac{rb}{R} + \frac{\sqrt{R^2 - a^2 - b^2}z_o}{R}\\
\end{bmatrix}.
\end{equation}
Note that $Q_{\vec x}$ is on the sphere such that $Q_{\vec x}\cdot Q_{\vec x} = R^2$. Equations (\ref{eqn4.8}) and
(\ref{eqn4.9}) represent transformations from global to local and local to global coordinates
respectively for each element.

The integration of edge departure points and quadrature points is performed on the sphere as length
preserving great circle arc, as is the intersection of swept regions and neighbouring elements.
In order to generate the quadrature points for the intersections however, these intersection 
regions must first be projected into the plane. We choose to use the plane centered on the intersecting
element in order to minimize errors due to the evaluation of points at large distances from the
tangent plane origin. There points are then projected back onto the sphere and integrated forwards in
time, where they are projected into the tangent plane of the target element $k$ and used to
evaluate its test functions as shown in the second term in the right hand side of (\ref{eqn1.11}).

\end{document}